\documentclass{amsart}
\usepackage[english]{babel}
\usepackage[cp1250]{inputenc}
\usepackage{amsmath,amssymb}
\usepackage{microtype}
\usepackage{float}
\usepackage{enumitem}
\usepackage{mathrsfs}

\usepackage{tikz}
\usetikzlibrary{matrix,arrows}
\usetikzlibrary{decorations.markings}

\newcommand\blfootnote[1]{%
  \begingroup
  \renewcommand\thefootnote{}\footnote{#1}%
  \addtocounter{footnote}{-1}%
  \endgroup
}

\newtheorem{thm}{Theorem}[section]
\newtheorem{prop}[thm]{Proposition}
\newtheorem{lem}[thm]{Lemma}
\newtheorem{cor}[thm]{Corollary}

\newtheorem*{claim*}{Claim}

\newtheorem{innercustomthm}{Theorem}
\newenvironment{customthm}[1]
  {\renewcommand\theinnercustomthm{#1}\innercustomthm}
  {\endinnercustomthm}

\theoremstyle{definition}
\newtheorem{definition}[thm]{Definition}
\newtheorem{remark}[thm]{Remark}
\newtheorem{example}[thm]{Example}

\newtheorem*{remark*}{Remark}

\def\FF{\mathbb F}

\def\NN{\mathbb N}

\DeclareMathOperator{\spinc}{Spin^c}

\DeclareMathOperator{\tb}{tb}
\DeclareMathOperator{\TB}{TB}
\DeclareMathOperator{\rot}{rot}
\DeclareMathOperator{\ROT}{ROT}
\DeclareMathOperator{\slk}{sl}

\DeclareMathOperator{\HFK}{HFK}

\begin{document}

\title[Non-loose negative torus knots]{Non-loose negative torus knots}

\author{Irena Matkovi\v{c}}
\address{Department of Mathematics, Uppsala University, Sweden}
\email{irena.matkovic@math.uu.se}

\begin{abstract}
We study Legendrian and transverse realizations of the negative torus knots $T_{(p,-q)}$ in all contact structures on the $3$-sphere. We give a complete classification of the strongly non-loose transverse realizations and the strongly non-loose Legendrian realizations with the Thurston-Bennequin invariant smaller than $-pq$. 

Additionally, we show that the strongly non-loose transverse realizations $T$ are classified by their non-zero invariants $\mathfrak T(T)$  in the minus version of the knot Floer homology. However, not all the elements of $\HFK^-(T_{(p,q)})$ can be  realized.

Along the way, we relate our Legendrian realizations to the tight contact structures on the Legendrian surgeries along them. Specifically, we realize all tight structures on the lens space $L(pq+1,p^2)$ as a single Legendrian surgery on a Legendrian $T_{(p,-q)}$, and we relate transverse realizations in overtwisted structures to the non-fillable tight structures on the large negative surgeries along the underlying knots.
\end{abstract}

%\subjclass[2010]{57R17}
%\subjclass[2020]{57K33, 57K18}

\maketitle

%%%%%%%%%%%%%%%%%%%%%%%%%%%%%%%%%%%%%%%%%%%%%%%%%%%%%%%%%%%%%%
\section{Introduction}

Legendrian and transverse knot theory in overtwisted structures was established as an interesting subject by the work of Etnyre \cite{E}: how a specific knot type is appearing in various overtwisted structures indicates significant differences among them. Furthermore, classification results bring to our attention contact counterparts of an important topological question asking which manifolds can be obtained by the surgery along a knot. On the other hand, the introduction of contact invariants in Heegaard Floer homology \cite{LOSSz} has raised the contact realization problem for its non-trivial classes. In the present paper, we intend to answer these questions for the first infinite family of knots -- negative torus knots -- building upon and confirming conjectures of Geiges and Onaran \cite{GO.L, GO.T}, and Lisca, Ozsv\'ath, Stipsicz and Szab\'o \cite{LOSSz}.

\subsubsection*{Knots in overtwisted structures}
The Legendrian and transverse knots in overtwisted contact manifolds are of two types depending on whether there is an overtwisted disk in the knot complement or not, we call them loose and non-loose, respectively. The level of non-looseness of a knot can be, as suggested by Baker and Onaran \cite{BO}, geometrically measured by the minimal number of intersections of the knot with any overtwisted disk in the manifold, called the depth, or by counting the number of stabilizations needed to loosen the knot, called the tension. 
 If the knot complement, additionally, has zero Giroux torsion (namely, less than $\pi$), we say that such a knot is strongly non-loose. 

Another subtlety of the Legendrian knots in overtwisted structures is that their classification up to Legendrian isotopy does not necessarily coincide with the classification up to contactomorphism, and the same holds for the transverse knots. The majority of the rare classification results in the literature limit themselves on the understanding of the contactomorphism type, usually called the coarse classification, and so we will do in the present paper. In the coarse setting the complete classification has been obtained for the loose knots, due to Etnyre \cite{E} classified by the (so-called) classical invariants, and for the unknot by Eliashberg and Fraser \cite{EF}. But, even in these simplest examples the classification does not go over to the isotopy level as shown by Vogel \cite{V}; in fact, to achieve this we would need some additional conditions on the position of overtwisted disks as in Dymara \cite{D} and Cavallo \cite{C}.

\subsubsection*{Torus knots}
What makes the study of torus knots accessible (also in the contact setting), is the fact that the knot complement is Seifert fibered; see Section \ref{SecSF} for details. In particular, this makes an array of  arguments, well-established in the case of closed Seifert manifolds, applicable to the study of non-loose representatives of torus knots. Also, for the study of torus knots we have an advantage of the classification being settled in the standard contact structure, owing to the work of Etnyre and Honda \cite{EH}.

Building on the above, we obtain here a classification of non-loose Legendrian and transverse negative torus knots, giving an explicit description for a representative of every equivalence class. Precedingly, only very limited cases have been studied, such as the case of the left-handed trefoil in the paper of Geiges and Onaran \cite{GO.T}.

\begin{thm}\label{Class}
Up to the boundary parallel Giroux torsion, every Legendrian negative torus knot $T_{(p,-q)}$ with tight complement and the Thurston-Bennequin invariant smaller than $\TB=-pq$ can be represented by some stabilization of a knot $L$ of Figure \ref{fig:LSdgm} which satisfies certain explicit conditions on the rotation numbers of the surgery curves (stated in Corollary \ref{cor:class}).

The transverse negative torus knots with the zero Giroux torsion complement are exactly transverse push-offs of those Legendrian torus knots whose every negative stabilization satisfies the above conditions.
\end{thm}

The precise classifications will be carried out in Theorem \ref{theorem:class} and Corollary \ref{cor:class} for the Legendrian negative torus knots, and in Corollary \ref{cor:transverse} for the transverse ones.  In particular (as we will observe in Remark \ref{rmk:d3}), all these non-loose realizations of $T_{(p,-q)}$ appear in the overtwisted structures whose $3$-dimensional invariant is positive, even, and at most $(p-1)(q-1)$; up to stabilizations there is at most one transversely non-loose realization of $T_{(p,-q)}$ in each structure. Note that since $T_{(p,-q)}$ is a fibered knot and it is not strongly quasi-positive, the open book with the binding $T_{(p,-q)}$ supports an overtwisted structure; and, this is the structure whose $3$-dimensional invariant equals $d_3=(p-1)(q-1)$.

\subsubsection*{Knot Floer homology}
Throughout the paper we will assume some basic knowledge of the knot Floer homology (as defined in \cite{OSz,R}); in particular, we will use the minus knot Floer homology of the torus knots. Furthermore, we recall that Lisca, Ozsv\'ath, Stipsicz and Szab\'o in \cite{LOSSz} defined the Legendrian invariant $\mathfrak{L}(L)$ of the null-homologous Legendrian knot $L\subset (Y,\xi)$, lying in the $\HFK^-(-Y,L)$ (so, the knot Floer homology of the mirror image in the case of the ambient manifold being the $3$-sphere). The invariant is known to be invariant under negative stabilizations, hence giving rise to an invariant of transverse knots, and is multiplied by $U$ by every positive stabilization. Furthermore, when the ambient contact manifold has non-vanishing contact invariant (so, for the sphere when we are in the standard contact structure) the invariant is non-vanishing for every Legendrian knot, and so, in particular, it has infinite $U$-order. In the overtwisted ambient, however, the invariant always has finite $U$-order, and it might vanish; in particular, it vanishes for all loose knots (and more, whenever there is Giroux $2\pi$-torsion, as observed by Stipsicz and V\'ertesi \cite{SV}).  In \cite{BO}, Baker and Onaran as another measure of the non-looseness suggest the order, defined as the sum of the $U$-torsion orders of Legendrian invariants for the knot and its orientation reverse. Finally, we recall from the work of Ozsv\'ath and Stipsicz \cite{OS} that the bigrading of the knot Floer homology group in which the invariant lies can be computed from the classical invariants of the Legendrian knot $L$ as \[2A(\mathfrak{L}(L))=\tb(L)-\rot(L)+1 \text{ and } M(\mathfrak{L}(L))=-d_3(\xi)+2A(\mathfrak{L}(L)).\]
where $A$ is the Alexander grading and $M$ is the Maslov grading.

One of the motivating questions for the present study has been to find out whether all the Heegaard Floer classes admit contact interpretation. 

We will prove (in Theorem \ref{theorem:invL}) that the strongly non-loose transverse $T_{(p,-q)}$ are classified by the non-vanishing invariants; in fact, they will be classified already by their Alexander grading, and hence, simply, by the self-linking number. 

\begin{thm}\label{transverse}
A transverse realization $T$ of the knot $T_{(p,-q)}$ in an overtwisted $S^3$ is strongly non-loose, if and only if it has non-zero $\mathfrak T(T)\in \HFK^-(T_{(p,q)})$. Two such knots are equivalent, if and only if they share the same invariant.
\end{thm}

For the special family of knots $T_{(2,-2n+1)}$ (see Example \ref{example}), our theorem confirms the conjectured non-vanishing of the Legendrian invariants for Legendrian knots presented in \cite[Section 6]{LOSSz}, by Lisca, Ozsv\'ath, Stipsicz and Szab\'o. In this case, the Legendrian invariants provide generators for the torsion part of the knot Floer homology. However, we will observe (by examples) that generally not all the $U$-torsion elements can be presented by the Legendrian invariants, not even in the overtwisted structure (equivalently, in the torsion summand $ \left(\FF[U]/(U^n)\right)_{(A,M)}$ with $d_3=2A-M$) in which there exists a transversely non-loose realization.

\subsubsection*{Tight contact structures on small Seifert manifolds}
Our understanding of the Legendrian torus knots with tight complements is built on the embedding into and comparison with the tight contact structures on the closed manifolds obtained by the contact surgery along these knots. Specifically, we will make use of the classification of tight and fillable contact structures on small Seifert fibered $L$-spaces $M(-1;r_1,r_2,r_3)$, given in \cite{M.t} and \cite{M.f} by the author. 

First, we use these classification results in order to put bounds on and to distinguish between tight contact structures on the knot complements. But, eventually they lead us to some intriguing relations between the tight contact structures on the knot complement and the tight contact structures on the surgeries along that knot.

\begin{prop}\label{Compare}
All the tight contact structures on the very negative integral surgeries along a negative torus knot arise by Legendrian surgery:  the fillable structures from the Legendrian realizations in the standard tight contact structure, and the tight non-fillable structures from the non-loose realizations in overtwisted structures. They can be counted in terms of the number of transverse realizations with zero Giroux torsion.
\end{prop}

Proposition \ref{Compare} refers to the statement of Proposition \ref{prop:number} and the observations made in the paragraphs above it, based on the classification result of Theorem \ref{theorem:class}.

Additionally, we observe the following result about Legendrian lens space surgeries, completing the work of Geiges and Onaran from \cite{GO.L}.

\begin{thm}\label{theoremL}
For any pair of coprime positive integers $p<q$, every tight contact structure on the lens space $L(pq+1,p^2)$ can be obtained by a single Legendrian surgery along some Legendrian realization of the negative torus knot $T_{(p,-q)}$ in some contact structure on $S^3$.
\end{thm}

\subsection*{Overview}
Section \ref{SecSF} elaborates on the Seifert fibered structure of the torus knot complements, and presents a way to equip them with contact structures. In Section \ref{SecL}, we prove Theorem \ref{theoremL} concerning tight structures on lens spaces. In Section \ref{SecHF}, we obtain classifications of non-loose negative torus knots (Theorem \ref{Class}), together with the non-vanishing of the knot Floer invariants (Theorem \ref{transverse}) and relations to the contact structures on the surgeries along these knots (Proposition \ref{Compare}).

\subsection*{Acknowledgement}
I am grateful for many helpful and motivating conversations I had with Sa\v{s}o Strle, Alberto Cavallo, and Andr\'as Stipsicz. My sincere thanks to Hyunki Min and John Etnyre for spotting a mistake in an early version. This research was supported by the European Research Council (ERC) under the European Union's Horizon 2020 research and innovation programme (grant agreement No 674978).

%%%%%%%%%%%%%%%%%%%%%%%%%%%%%%%%%%%%%%%%%%%%%%%%%%%%%%%%%%%%%%
\section{Seifert fibration of the knot complement}\label{SecSF}

Let $p$ and $q$ be positive integers such that $p<q$ and $\gcd(p,q)=1$, and write $T_{(p,-q)}$ for the negative $(p,-q)$-torus knot. 

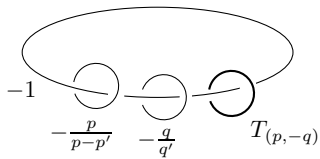
\begin{figure}[h]%----------------------------------------------------------------------------
\centering
\begin{tikzpicture}
\draw[white] (7.5,-0.5) circle (0pt) node[black,left,scale=0.9]{$-1$};
\draw (9,0) ellipse (1.8cm and 0.6cm);
\draw [fill=white,white] (8.35,-0.45) rectangle (8.55,-0.65);
\draw [fill=white,white] (9.3,-0.6) rectangle (9.45,-0.4);
\draw [fill=white,white] (10.1,-0.35) rectangle (10.3,-0.55);
\draw (7.9,-.55) arc (-160:170:0.3);
\draw[white] (8,-0.8) circle (0pt) node[below,black,scale=0.9]{$-\frac{p}{p-p'}$};
\draw (8.8,-.7) arc (-160:160:0.3);
\draw[white] (9,-.9) circle (0pt) node[below,black,scale=0.9]{$-\frac{q}{q'}$};
\draw[thick] (9.7,-.65) arc (-160:170:0.3);
\draw[white] (10.7,-.8) circle (0pt) node[below,black,scale=0.9]{$T_{(p,-q)}$};
\end{tikzpicture}
\caption{Torus knot $T_{(p,-q)}$.}
\label{fig:SFS}
\end{figure}

It is well-known that the complement of a torus knot is Seifert fibered. Concretely, the complement of the knot $T_{(p,-q)}$ is Seifert fibered over the disk with two singular fibers whose Seifert invariants are $-\frac{p'}{p}$ and $-\frac{q'}{(-q)}$ for $p',q'$ such that $pq'-qp'=1$. Since $-1<-\frac{p'}{p}<0$ and $\frac{q'}{q}>0$, we can renormalize the invariants as in Figure \ref{fig:SFS}; denote this manifold as $M(D^2; \frac{p-p'}{p},\frac{q'}{q})$.

%--------------------------------------------------------------------------------

If we write $q=np-k$ for $n\geq 2$ and $k<p,\ \gcd(p,k)=1$, the Seifert invariants equal $M(D^2;\frac{p-C}{p},\frac{Cn-D}{np-k})$ where $C,D$ are positive integers satisfying $Ck=Dp+1$. Furthermore, the negative continued fraction expansions of the two invariants are related as follows; here, we use the convention
 \[ [c_0,\ldots,c_m]=c_0-\frac{1}{\ddots-\frac{1}{c_m}}. \]

\begin{lem}\label{lem:cfe}
With positive integers $p,n,k,C,D$ chosen as above, the negative continued fraction expansions of the two Seifert invariants equal
\[ \frac{p}{p-C}=[a_1^0,\dots,a_1^s] \text{ and } \frac{np-k}{Cn-D}=[a_2^0,\dots,a_2^t,n] \]
for some $a_i^j\geq2$ satisfying $[a_1^0,\dots,a_1^s]^{-1}+[a_2^0,\dots,a_2^t]^{-1}=1$.
\end{lem}

\begin{proof}
If we write out \[ \frac{p}{p-C}=[a_1^0,\dots,a_1^{m_1}] \text{ and } \frac{np-k}{Cn-D}=[a_2^0,\dots,a_2^{m_2}], \]
the chain of unknots with coefficients $(-a_1^{m_1},\dots,-a_1^0,-1,-a_2^0,\dots,-a_2^{m_2})$ gives a surgery description of the ambient manifold, which is $S^3$.

On the other hand, since $\frac{p-C}{p}+\frac{Cn-D}{np-k}>1$, there are truncated continued fractions such that $[a_1^0,\dots,a_1^s]^{-1}+[a_2^0,\dots,a_2^t]^{-1}=1$ for $s\leq m_1$ and $t\leq m_2$, and the truncated chains of unknots join into a chain with surgery coefficients $(-a_1^s,\dots,-a_1^0,-1,-a_2^0,\dots,-a_2^t)$ which corresponds to $S^1\times S^2$ (for more details, see \cite[Lemma 3.1]{M.f}). 

Now, the only way to get $S^3$ from $S^1\times S^2$ by lengthening the chain is by adding a single unknot at one of the two ends. Since the numerator of the second fraction is larger, the coefficient is added to the second continued fraction; so, $m_1=s$ and $m_2=t+1$.  Finally, that the added coefficient equals $n$ can be seen from the comparison to the continued fraction expansion of $\frac{np^2-kp+1}{p^2}$ which starts in $n$ (see the proof of Theorem \ref{theorem:L}).
\end{proof}

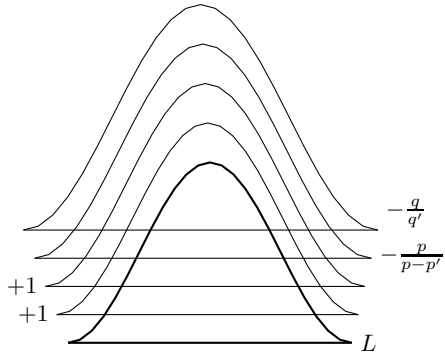
\begin{figure}[h]%----------------------------------------------------------------------------
\centering
\begin{tikzpicture}
\begin{scope}[scale=1.5]
\begin{scope}[shift={(0,0)}]
\draw [domain=-0.5*pi:1.5*pi, scale=0.5] plot (\x, {2*sin(\x r)});
\draw [scale=0.5] (-0.5*pi,-2)--(1.5*pi,-2);
\draw [white,scale=0.5] (1.5*pi,-1.7) circle (0.1pt) node[black,right,scale=0.9]{$-\frac{q}{q'}$};
\end{scope}
\begin{scope}[shift={(.06,-.3)},scale=0.95]
\draw [domain=-0.5*pi:1.5*pi, scale=0.5] plot (\x, {2*sin(\x r)});
\draw [scale=0.5] (-0.5*pi,-2)--(1.5*pi,-2) node[right,scale=0.9]{$-\frac{p}{p-p'}$};
\end{scope}
\begin{scope}[shift={(.12,-.6)},scale=0.9]
\draw [domain=-0.5*pi:1.5*pi, scale=0.5] plot (\x, {2*sin(\x r)});
\draw [scale=0.5] (-0.5*pi,-2)--(1.5*pi,-2);
\draw [scale=0.5] (-0.5*pi,-2) circle (0.1pt) node[left,scale=0.9]{$+1$};
\end{scope}
\begin{scope}[shift={(.18,-.9)},scale=0.85]
\draw [domain=-0.5*pi:1.5*pi, scale=0.5] plot (\x, {2*sin(\x r)});
\draw [scale=0.5] (-0.5*pi,-2)--(1.5*pi,-2);
\draw [scale=0.5] (-0.5*pi,-2) circle (0.1pt)node[left,scale=0.9]{$+1$};
\end{scope}
\begin{scope}[shift={(.24,-1.2)},scale=0.8]
\draw [thick,domain=-0.5*pi:1.5*pi, scale=0.5] plot (\x, {2*sin(\x r)});
\draw [thick,scale=0.5] (-0.5*pi,-2)--(1.5*pi,-2)node[right,scale=0.9]{$L$};
\end{scope}
\end{scope}
\end{tikzpicture}
\caption{Legendrian realizations of the torus knot $T_{(p,-q)}$.}
\label{fig:LSdgm}
\end{figure}

Some contact structures on the above (bounded) Seifert manifolds can be described by contact surgery diagrams of Figure \ref{fig:LSdgm}. These diagrams first appeared in the work of Lisca and Stipsicz \cite{LS.III}, and have since been extensively used in understanding tight structures on Seifert fibered spaces \cite{LS, GLS, M.t, M.f}, as well as for providing examples of non-loose knots \cite{LOSSz, GO.L, GO.T}.

%---------------------------------------------------------------------------------

Recall that such a diagram gives a family of contact structures, whose elements can be specified by replacing each rational contact surgery with a Legendrian surgery along a chain of unknots whose Thurston-Bennequin invariants are determined by the continued fraction expansion
as \[\tb_i^0=-a_i^0\  \text{ and}\  \tb_i^j=-a_i^j+1 \text{ for } j>0, \]
and rotation numbers are chosen arbitrarily in \[\rot_i^j\in\{\tb_i^j+1,\tb_i^j+3,\dots,-\tb_i^j-1\}.\]

So, a single Legendrian representation of $T_{(p,-q)}$ from Figure \ref{fig:LSdgm} is specified by the array of rotation numbers for the unknots supporting Legendrian surgery.

Let us recall some terminology from \cite{M.f}.
\begin{definition}
We say that a Legendrian unknot is fully positive if all its stabilizations are positive, that is $\rot = -(\tb +1)$. Analogously, a Legendrian unknot is fully negative if all its stabilizations are negative, that is $\rot = \tb +1$.

Additionally, a contact surgery along four $(-1)$-linked unstabilized Legendrian unknots (as in Figure \ref{fig:LSdgm}) with $(+1)$-surgery performed along two of them and the negative inverses of the rational surgery coefficients of the other two adding to one, is called a balanced link when turned into a Legendrian surgery along the chains of Legendrian unknots, all the unknots of one chain are fully positive and all the unknots of the other chain are fully negative.
\end{definition}

\begin{prop}\label{prop:tightC}
For any choice of rotation numbers, Figure \ref{fig:LSdgm} presents Legendrian torus knot $T_{(p,-q)}$ with tight complement in some contact structure on $S^3$. Moreover, the ambient contact structure on $S^3$ is tight if and only if the contact surgery presentation contains a balanced link.
\end{prop}

\begin{proof}
For tightness of the knot complement we use the standard cancellation argument: since the $(-1)$-surgery along $L$ results in a tight contact manifold, the complement of $L$ cannot be overtwisted.

Knowing that the only tight structure on $S^3$ is also Stein fillable, the question whether the contact structure on the ambient $S^3$, given by the surgery diagram of Figure \ref{fig:LSdgm}, is tight or overtwisted, is in greater generality answered in \cite[Theorem 1.1]{M.f}. It is equivalent to the surgery presentation containing a balanced sublink, which in our case is fulfilled by either
\[\rot_1^j=\tb_1^j+1 \text{ for all } j\in\{0,\dots,s\} \text{ and } \rot_2^j=-\tb_2^j-1 \text{ for }j\in\{0,\dots,t\},\]
or
\[\rot_1^j=-\tb_1^j-1 \text{ for all } j\in\{0,\dots,s\} \text{ and } \rot_2^j=\tb_2^j+1 \text{ for }j\in\{0,\dots,t\}.\]
\end{proof}

%%%%%%%%%%%%%%%%%%%%%%%%%%%%%%%%%%%%%%%%%%%%%%%%%%%%%%%%%%%%%%
\section{Legendrian knots with  $tb = -pq$ and tight $L(pq+1,p^2)$} \label{SecL}

\begin{lem}\label{lem:tb}
In any contact surgery presentation of Figure \ref{fig:LSdgm}, the Legendrian realization $L$ of the torus knot $T_{(p,-q)}$ has the Thurston-Bennequin invariant equal $\tb(L)=-pq$.
\end{lem}

\begin{proof}
We use the formula from \cite[Lemma 6.6]{LOSSz}: \[ \tb(L)= \tb_0 + \frac{\det(Q(0,0,0,x,y))}{\det(Q(0,0,x,y))} \]
where $\tb_0$ is the Thurston-Bennequin invariant of the knot before surgery, and $Q$ with $x=-1-\frac{q}{q'}, y=-1-\frac{p}{p-p'}$ is the intersection matrix of the (smooth) surgery diagram of Figure \ref{fig:LSdgm}. Then:
\[ \tb(L) = -1 + \frac{-4-3x-3y-2xy}{-3-2x-2y-xy} = \frac{-pq}{pq'-p'q} =-pq. \]
\end{proof}

We know (since Moser \cite{Mo}) that all torus knots are lens space knots; explicitly, the $-(pq\pm1)$-surgery along the negative torus knot $T_{(p,-q)}$ results in the lens space $L(pq\pm1,p^2)$. 

In \cite{GO.L}, Geiges and Onaran studied Legendrian lens space surgeries, culminating in a presentation of every tight contact structure on $L(np^2-p+1,p^2)$ as a Legendrian surgery on some Legendrian realization of $T_{(p,-(np-1))}$. 
We generalize their result to every negative torus knot, completely confirming the conjecture stated in \cite[Remark 1.2 (3)]{GO.L}.

\begin{customthm}{\ref{theoremL}}\label{theorem:L}
For any pair of coprime positive integers $p<q$, every tight contact structure on the lens space $L(pq+1,p^2)$ can be obtained by a single Legendrian surgery along some Legendrian realization of the negative torus knot $T_{(p,-q)}$ in some contact structure on $S^3$.
\end{customthm}

\begin{proof}
Since in Lemma \ref{lem:tb} we have computed the Thurston-Bennequin invariant of any Legendrian realization $L$ of $T_{(p,-q)}$ from Figure \ref{fig:LSdgm} to be $-pq$, we know that Legendrian surgery along any such $L$ results in some contact $L(pq+1,p^2)$. In fact, as noticed already in the proof of Proposition \ref{prop:tightC}, the resulting contact structure is tight. We need to show that by varying Legendrian realization $L$ -- that is, by choosing different rotation numbers for the surgery curves -- we reach all tight contact structures on $L(pq+1,p^2)$.

As specified by Honda \cite{Ho.I}, the tight structures on a lens space $L(u,v)$ for relatively prime $u$ and $v$ correspond to the choices of rotation numbers on the chain of Legendrian unknots whose Thurston-Bennequin invariants are determined by the negative continued fraction expansion of $\frac{u}{v}$. In fact, the tight structures are distinguished already by their induced $\spinc$-structures.

Geiges and Onaran (see \cite[Theorem 1.1]{GO.L}) proved the theorem in the case $q=np-1$ by finding as many different rotation numbers for Legendrian $T_{(p,-q)}$ with $\tb=-pq$ as there are tight contact structures on  $L(pq+1,p^2)$. In contrast, we obtain a direct comparison of the two surgery presentations: the standard one as a Legendrian surgery along the chain of unknots, and the one given by a single Legendrian surgery along a Legendrian $T_{(p,-q)}$ in some contact $S^3$. For the second, we consider Legendrian knots $L$ of Figure \ref{fig:LSdgm}.

\begin{figure}[h]%--------------------------------------------------------------------------
\centering
\begin{tikzpicture}\begin{scope}[scale=1.3]
\begin{scope}[shift={(-0.5,0)}]
\draw (1.3,-0.5) circle (0.9cm);
\draw [fill=white,white] (0.8,-1) rectangle (0.5,-1.1);
\draw (3.2,-0.5) circle (0.9cm);
\draw [fill=white,white] (3.7,-1) rectangle (4,-1.1);
\draw (1.3,-1.8) circle (0pt) node[above,scale=0.7]{$-a_1^0-1$};
\draw (.8,.35) circle (0pt) node[above,scale=0.7]{$(\rot_1^0)$};
\draw (3.2,-1.8) circle (0pt) node[above,scale=0.7]{$-a_2^0-1$};
\draw (3.7,.35) circle (0pt) node[above,scale=0.7]{$(\rot_2^0)$};
\draw (2.9,-0.5) circle (0.5cm);
\draw (3.25,-1) circle (0pt) node[right,scale=0.7]{$0$};
\draw (3.25,-0.1) circle (0pt) node[right,scale=0.7]{$(0)$};
\draw [fill=white] (2.05,-0.62) rectangle (2.5,-0.3) node[below left,scale=0.7]{$-1$};

\end{scope}
\begin{scope}[shift={(3.72,-1)}]
\draw (1,0) ++(40:.33 and .2) arc (40:-220:.33 and 0.2)
          (1,0) ++(130:.33 and .2) arc (130:50:.33 and 0.2);
\draw (0,0) ++(40:.33 and .2) arc (40:-220:.33 and 0.2)
          (0,0) ++(110:.33 and .2) arc (110:50:.33 and 0.2);
\draw (0.5,0) ++(-20:.3) arc (-20:200:.3)
          (0.5,0) ++(-145:.3) arc (-145:-35:.3)
          (1.5,0) ++(-145:.3) arc (-145:200:.3); 
\draw[white] (0,-.21) circle (0pt) node[below,black,scale=0.7]{$-a_2^1$};
\draw [fill=white,white] (.4,-.4) rectangle (.6,.4); \draw[white] (.5,.1) circle (0pt) node[below,black,scale=.4]{$\cdots$};
\draw[white] (1,-.21) circle (0pt) node[below,black,scale=0.7]{$-a_2^t$};
\draw[white] (1.5,-.3) circle (0pt) node[below,black,scale=0.7]{$-n$};
\draw[white] (0.1,.21) circle (0pt) node[above,black,scale=0.6]{$(\rot_2^1)$};
\draw[white] (1,.21) circle (0pt) node[above,black,scale=0.6]{$(\rot_2^t)$};
\draw[white] (1.5,.3) circle (0pt) node[above,black,scale=0.6]{$(\rot_2^{t+1})$};
\end{scope}

\begin{scope}[shift={(-1.72,-1)}]
\draw (0,0) ++(-20:.3) arc (-20:325:.3);
\draw (.5,0) ++(40:.33 and .2) arc (40:-220:.33 and 0.2)
          (.5,0) ++(130:.33 and .2) arc (130:50:.33 and 0.2); 
\draw (1,0) ++(-20:.3) arc (-20:200:.3)
          (1,0) ++(-145:.3) arc (-145:-35:.3);
\draw (1.5,0) ++(40:.33 and .2) arc (40:-220:.33 and 0.2)
          (1.5,0) ++(130:.33 and .2) arc (130:70:.33 and 0.2);
\draw [fill=white,white] (.4,-.3) rectangle (.6,.3); \draw[white] (.5,.1) circle (0pt) node[below,black,scale=.4]{$\cdots$};
\draw[white] (0,-.3) circle (0pt) node[below,black,scale=0.7]{$-a_1^s$};
\draw[white] (1,-.3) circle (0pt) node[below,black,scale=0.7]{$-a_1^2$};
\draw[white] (1.5,-.21) circle (0pt) node[below,black,scale=0.7]{$-a_1^1$};

\draw[white] (0,.3) circle (0pt) node[above,black,scale=0.6]{$(\rot_1^s)$};
\draw[white] (1,.3) circle (0pt) node[above,black,scale=0.6]{$(\rot_1^2)$};
\draw[white] (1.41,.21) circle (0pt) node[above,black,scale=0.6]{$(\rot_1^1)$};
\end{scope}

%-------------------------------------------------------------
\begin{scope}[shift={(.45,0)}]
\begin{scope}[shift={(-0.5,-3)}]
\begin{scope}[decoration={markings,mark=at position .011 with {\arrow{>}}}]
\draw [postaction={decorate}] (.55,-0.5) circle (0.9cm);
\end{scope}
\draw [fill=white,white] (-.3,-1) rectangle (0.2,-1.1);
\draw (.85,-1.7) circle (0pt) node[above,scale=0.7]{$-a_1^0$};
\begin{scope}[decoration={markings,mark=at position 0.5 with {\arrow{<}}}]
\draw [postaction={decorate}] (3.2,-0.5) circle (0.9cm);
\end{scope}
\draw [fill=white,white] (3.7,-1) rectangle (4,-1.1);
\draw (3.5,-1.7) circle (0pt) node[above,scale=0.7]{$-a_2^0$};
\draw (-.1,.35) circle (0pt) node[above,scale=0.7]{$(\rot_1^0+1)$};
\draw (3.85,.35) circle (0pt) node[above,scale=0.7]{$(\rot_2^0+1)$};
\draw (1.88,.1) circle (0pt) node[above,scale=0.7]{$(0)$};
\draw [fill=white,white] (1.3,-0.9) rectangle (2.5,-0.7);
\draw (1.88,-.5) ++(40:.66 and .4) arc (40:-220:.66 and 0.4)
          (1.88,-.5) ++(130:.66 and .4) arc (130:50:.66 and 0.4);
\draw (1.88,-1.3) circle (0pt) node[above,scale=0.7]{$0$};

\end{scope}
\begin{scope}[shift={(3.72,-4)}]
\draw (1,0) ++(40:.33 and .2) arc (40:-220:.33 and 0.2)
          (1,0) ++(130:.33 and .2) arc (130:50:.33 and 0.2);
\draw (0,0) ++(40:.33 and .2) arc (40:-220:.33 and 0.2)
          (0,0) ++(110:.33 and .2) arc (110:50:.33 and 0.2);
\draw (0.5,0) ++(-20:.3) arc (-20:200:.3)
          (0.5,0) ++(-145:.3) arc (-145:-35:.3)
          (1.5,0) ++(-145:.3) arc (-145:200:.3); 
\draw[white] (0,-.21) circle (0pt) node[below,black,scale=0.7]{$-a_2^1$};
\draw [fill=white,white] (.4,-.4) rectangle (.6,.4); \draw[white] (.5,.1) circle (0pt) node[below,black,scale=.4]{$\cdots$};
\draw[white] (1,-.21) circle (0pt) node[below,black,scale=0.7]{$-a_2^t$};
\draw[white] (1.5,-.3) circle (0pt) node[below,black,scale=0.7]{$-n$};
\end{scope}

\begin{scope}[shift={(-2.48,-4)}]
\draw (0,0) ++(-20:.3) arc (-20:325:.3);
\draw (.5,0) ++(40:.33 and .2) arc (40:-220:.33 and 0.2)
          (.5,0) ++(130:.33 and .2) arc (130:50:.33 and 0.2); 
\draw (1,0) ++(-20:.3) arc (-20:200:.3)
          (1,0) ++(-145:.3) arc (-145:-35:.3);
\draw (1.5,0) ++(40:.33 and .2) arc (40:-220:.33 and 0.2)
          (1.5,0) ++(130:.33 and .2) arc (130:70:.33 and 0.2);
\draw [fill=white,white] (.4,-.3) rectangle (.6,.3); \draw[white] (.5,.1) circle (0pt) node[below,black,scale=.4]{$\cdots$};
\draw[white] (0,-.3) circle (0pt) node[below,black,scale=0.7]{$-a_1^s$};
\draw[white] (1,-.3) circle (0pt) node[below,black,scale=0.7]{$-a_1^2$};
\draw[white] (1.5,-.21) circle (0pt) node[below,black,scale=0.7]{$-a_1^1$};
\end{scope}
\end{scope}

%-------------------------------------------------------------
\begin{scope}[shift={(-.77,-3)}]
\begin{scope}[shift={(-0.5,-3)}]
\draw (3.2,-0.5) circle (0.9cm);
\draw [fill=white,white] (1.7,-1) rectangle (4,-1.1);
\draw (3.2,-1.8) circle (0pt) node[above,scale=0.7]{$-a_1^0-a_2^0$};
\draw (3.2,.45) circle (0pt) node[above,scale=0.7]{$(\rot_1^0-\rot_2^0)$};

\end{scope}
\begin{scope}[shift={(3.72,-4)}]
\draw (1,0) ++(40:.33 and .2) arc (40:-220:.33 and 0.2)
          (1,0) ++(130:.33 and .2) arc (130:50:.33 and 0.2);
\draw (0,0) ++(40:.33 and .2) arc (40:-220:.33 and 0.2)
          (0,0) ++(110:.33 and .2) arc (110:50:.33 and 0.2);
\draw (0.5,0) ++(-20:.3) arc (-20:200:.3)
          (0.5,0) ++(-145:.3) arc (-145:-35:.3)
          (1.5,0) ++(-145:.3) arc (-145:200:.3); 
\draw[white] (0,-.21) circle (0pt) node[below,black,scale=0.7]{$-a_2^1$};
\draw [fill=white,white] (.4,-.4) rectangle (.6,.4); \draw[white] (.5,.1) circle (0pt) node[below,black,scale=.4]{$\cdots$};
\draw[white] (1,-.21) circle (0pt) node[below,black,scale=0.7]{$-a_2^t$};
\draw[white] (1.5,-.3) circle (0pt) node[below,black,scale=0.7]{$-n$};
\end{scope}

\begin{scope}[shift={(.18,-4)}]
\draw (0,0) ++(-20:.3) arc (-20:325:.3);
\draw (.5,0) ++(40:.33 and .2) arc (40:-220:.33 and 0.2)
          (.5,0) ++(130:.33 and .2) arc (130:50:.33 and 0.2); 
\draw (1,0) ++(-20:.3) arc (-20:200:.3)
          (1,0) ++(-145:.3) arc (-145:-35:.3);
\draw (1.5,0) ++(40:.33 and .2) arc (40:-220:.33 and 0.2)
          (1.5,0) ++(130:.33 and .2) arc (130:70:.33 and 0.2);
\draw [fill=white,white] (.4,-.3) rectangle (.6,.3); \draw[white] (.5,.1) circle (0pt) node[below,black,scale=.4]{$\cdots$};
\draw[white] (0,-.3) circle (0pt) node[below,black,scale=0.7]{$-a_1^s$};
\draw[white] (1,-.3) circle (0pt) node[below,black,scale=0.7]{$-a_1^2$};
\draw[white] (1.5,-.21) circle (0pt) node[below,black,scale=0.7]{$-a_1^1$};
\end{scope}
\end{scope}

\end{scope}
\end{tikzpicture}
\caption{Kirby diagrams for $L(np^2-kp+1,p^2)$ with the $\spinc$-structure.}
\label{fig:L}
\end{figure}
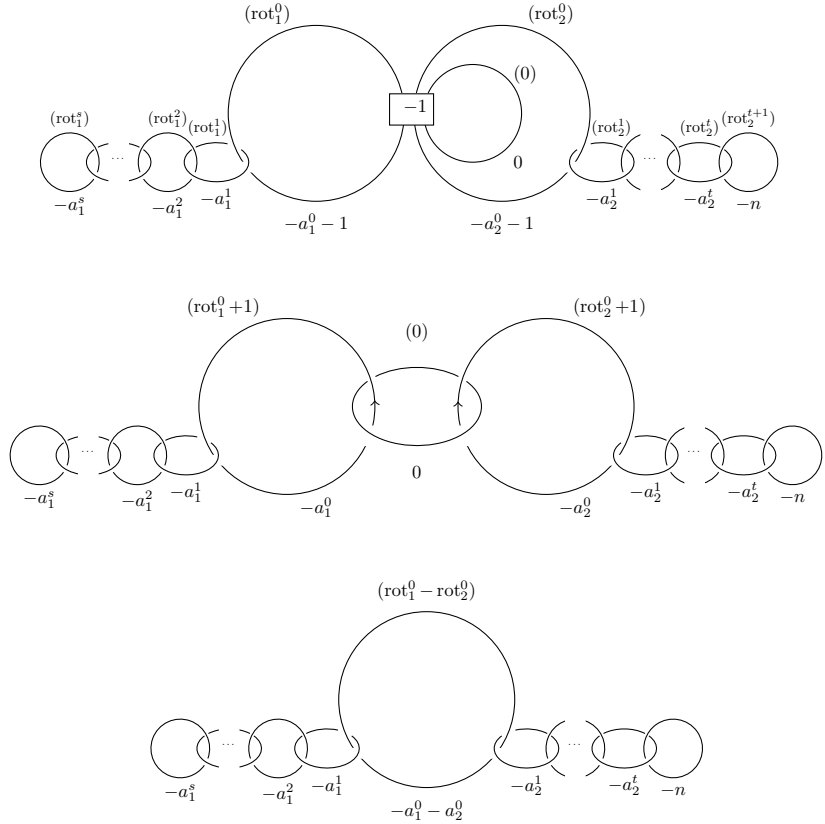%--------------------------------------------------------------------------

Together with the Legendrian surgery along $L$, the surgery diagram of Figure \ref{fig:LSdgm} smoothly describes $L(pq+1,p^2)$ and it looks like the first diagram of Figure \ref{fig:L}. In Figure \ref{fig:L}, we use $q=np-k$ as in Lemma \ref{lem:cfe} and we keep track of the $\spinc$-structure, induced by the chosen contact structure, by writing its evaluations on the homology generators (in the parenthesis above the corresponding knots). 

We get from the first to the second diagram by blowing up the $(-1)$-linking point followed by a blow-down of the $(+1)$-framed meridian of the thus-added curve. Then, from the second to the third diagram we get by sliding the $(-a_1^0)$-framed unknot over the reversely oriented $(-a_2^0)$-framed unknot, thus unlinking the $(-a_1^0)$-framed unknot from the $0$-framed unknot, and a consecutive cancellation of the $(-a_2^0)$-framed unknot with its $0$-framed meridian.

Now, since $[a_1^0,\dots,a_1^s]^{-1}+[a_2^0,\dots,a_2^t]^{-1}=1$, we know that one of the two coefficients $a_1^0$ or $a_2^0$ equal $2$. Hence, we can reach all possible rotation numbers on the $(-a_1^0-a_2^0)$-framed unknot in the chain by differences of the initial rotation numbers. Indeed, without loss of generality let us write $a_1^0=2$ and $a_2^0=m$, then the possible rotation numbers for the $(-2-m)$-framed unknot are $-m, -m+2,\dots,m-2, m$, and we can get them by choosing, for example, the pairs of rotation numbers $(-1,m-1), (-1,m-3),\dots,(-1,-m+1),(1,m+1)$ for the $(-3)$- and $(-m-1)$-framed unknots in the initial diagram. The choice of all the other rotation numbers can be taken equal in both diagrams. So, we have realized all possible choices of the rotation numbers on the chain by possible choices of rotation numbers in our surgery diagrams, and by that, all possible tight contact structures on $L(pq+1,p^2)$ by a Legendrian surgery along $L$ in some contact structure on $S^3$.
\end{proof}

%%%%%%%%%%%%%%%%%%%%%%%%%%%%%%%%%%%%%%%%%%%%%%%%%%%%%%%%%%%%%%
\section{Legendrian knots with $tb < -pq$ and knot Floer invariants} \label{SecHF}

Stabilizing a Legendrian knot does not change its knot type, hence any $\ell$-times stabilization of a knot $L$ from Figure \ref{fig:LSdgm} gives us a Legendrian $T_{(p,-q)}$ with $\tb=-pq-\ell$. Let us denote any $\ell$-times stabilization of any knot $L$ by $L^\ell$. As before a single Legendrian realization is specified by an array of rotation numbers for all the surgery curves in Figure \ref{fig:LSdgm}, and here additionally, by the rotation number of $L^\ell$. Be aware, however, that for $\ell\geq 1$ not all $L^\ell$ have tight complements.

Let us recall an equivalence relation on the set of Legendrian or transverse knots.

\begin{definition}
Two Legendrian or transverse knots $K_1$ and $K_2$ in a closed contact $3$-manifold $(M,\xi)$ are said to be (coarse) equivalent if there exists a contactomorphism of $(M,\xi)$ which maps $K_1$ to $K_2$.
\end{definition}

When the self-contactomorphism group of $(M,\xi)$ is not contractible, the coarse equivalence is known (see for example \cite{V}) to be weaker than Legendrian isotopy. However, based on Kegel's \cite[Lemma 10.3]{K}, we know that in the case of nontrivial knots in the $3$-sphere, when we have no cosmetic surgeries, the equivalence type of a Legendrian knot is completely determined by the contactomorphisms of the knot complement, even when we look at all contact structures on $S^3$ simultaneously.

\begin{thm}[Coarse classification of Legendrian $T_{(p,-q)}$ with $\tb<-pq$] \label{theorem:class}
Let $K$ be a Legendrian realization of the knot $T_{(p,-q)}$ with $\tb(K)=-pq-\ell$ in some contact structure on $S^3$ . The complement of $K$ has Giroux torsion equal to zero, if and only if $K$ is equivalent to some $L^\ell$ along which the Legendrian surgery results in a tight contact manifold. 

Moreover, two Legendrian knots $K_1$ and $K_2$ are equivalent, if and only if the two contact manifolds obtained by Legendrian surgery along the knot $K_i$ and its Legendrian push-off for $i=1$ and $2$ are contact isotopic.
\end{thm}

\begin{proof}
We will first show that there is a one-to-one correspondence between the equivalence classes of Legendrian knots $T_{(p,-q)}$ with $\tb=-pq-\ell$ whose complement has zero Giroux torsion, and the tight contact structures on $M(-1;\frac{p-p'}{p},\frac{q'}{q},\frac{2}{2\ell+1})$ up to contact isotopy. This will prove that the knot complement has zero Giroux torsion, if and only if it is equivalent to some $L^\ell$ and the Legendrian surgery along it and its Legendrian push-off results in a tight manifold, and that two knots are equivalent, if and only if the results of surgery are isotopic. Then, we will observe that the contact manifold obtained by the Legendrian surgery along a Legendrian knot $L^\ell$ is tight, if and only if the contact manifold obtained by the Legendrian surgery along the knot and its Legendrian push-off is tight, thus completing the proof of the theorem.

If we write out the Seifert fibration of the $T_{(p,-q)}$ complement as $M(D^2;\frac{p-p'}{p},\frac{q'}{q})$, and choose a trivialization of its boundary torus by the meridian $-\partial D^2\times\{1\}$ and the longitude a parallel Seifert fibre $\{\ast\}\times S^1$ for $\ast\in\partial D^2$, then the slope of dividing curves on the boundary torus (once perturbed to be minimal convex) equals $s=\frac{1}{\tb-pq}+1$ (for details see \cite[Section 4]{GO.T}). Utilizing the idea of Ding, Li and Zhang \cite{DLZ}, we want to embed this knot complement in some small Seifert fibered manifold $M(-1;\frac{p-p'}{p},\frac{q'}{q},r_3)$ such that the boundary slope $s$ coincides with the slope $-1$ when measured in the standard basis for the neighborhood of the knot; by the formula in \cite[page 65]{DLZ}, the third Seifert constant then equals $r_3=\frac{2}{2\ell+1}=[\ell+1,2]^{-1}$.  Since there is a unique tight structure with boundary slope $-1$ on the solid torus, this immediately tells that on the $T_{(p,-q)}$ complement with boundary slope $s$ there are at least as many structures with zero Giroux torsion as on $M=M(-1;\frac{p-p'}{p},\frac{q'}{q},\frac{2}{2\ell+1})$. What is more, for both structures the maximal twisting number of the regular fiber (of the Seifert fibration) is equal to zero, hence, according to Lisca and Stipsicz \cite[Proposition 6.1]{LS.III} they can all be presented by surgery diagrams of Figure \ref{fig:LSdgm}. Examining the upper bound for the number of tight structures on $M$ in \cite[Section 5]{M.t}, we observe that all the overtwistedness as well as all the isotopies between different surgery presentations have been achieved in the complement of the singular fibers; hence, actually giving the upper bound for the number of structures with zero Giroux torsion on the knot complement. Here, the condition of the zero Giroux torsion comes from the fact that the boundary parallel Giroux $\pi$-torsion gives rise to an overtwisted disk once we do Legendrian surgery along the core knot, and hence in \cite{M.t} only appropriate structures on the circle bundle over the pair of pants were taken into consideration. This proves our first assertion about the correspondence with the tight contact structures on $M(-1;\frac{p-p'}{p},\frac{q'}{q},\frac{2}{2\ell+1})$, the relation is acomplished by Legendrian surgery along the knot and its push-off.

In order to prove the second assertion about the two surgeries being tight for the same Legendrian knots, we need to look closer at the isotopy and overtwistedness conditions in \cite[Section 5]{M.t}. The manifold obtained by the Legendrian surgery along a Legendrian knot $L^\ell$ is $M(-1;\frac{p-p'}{p},\frac{q'}{q},\frac{1}{\ell+1})$, and the manifold obtained by the Legendrian surgery along the knot and its Legendrian push-off is $M(-1;\frac{p-p'}{p},\frac{q'}{q},\frac{2}{2\ell+1})$. First, we notice that isotopies of type (I3) of \cite[Proposition 5.2]{M.t} and overtwisted structures of type (O2) of \cite[Proposition 5.1]{M.t} do not occur for either of considered manifolds, because the continued fraction expansions of the Seifert constants cannot fulfill the required equality. Hence, different surgery presentations of $M(-1;\frac{p-p'}{p},\frac{q'}{q},\frac{1}{\ell+1})$ are isotopic if they are related by either  (I2) or a sequence of (I1) changes of rotation numbers as in \cite[Proposition 5.2]{M.t}, while in the case of $M(-1;\frac{p-p'}{p},\frac{q'}{q},\frac{2}{2\ell+1})$ the relation (I2) does not give an isotopy. Anyway, the overtwisted structures in both cases are described by (O1) of \cite[Proposition 5.1]{M.t} or the structures which are related to a structure satisfying (O1) by a sequence of (I1) changes.
\end{proof}

In fact, we can write out explicitly when a Legendrian knot $L^\ell$ is non-loose and whether two $L^\ell$ are coarse equivalent. Recall that a Legendrian knot $L^\ell$ is specified by its rotation number and the rotation numbers for the surgery curves of Figure \ref{fig:LSdgm}. We will write the rotation numbers of the unknots forming the two singular fibers as $\rot_j^i$, and use $p_j^i$ and $n_j^i$ for the number of their positive and negative stabilizations respectively (of course $\rot_j^i=p_j^i-n_j^i$). Note though that we do not fix $j=1,2$ to a particular Seifert constant, but rather take $j=1$ for the singular fiber with the longer fully negative truncation and $j=2$ for the one with the longer fully positive truncation. (To be clear, here by longer we mean that the denominator of the inverse of the corresponding continued fraction is bigger; however, the inequality might not be strict for both positive and negative truncations, and then the one which is strict determines the order of singular fibers, or if neither is strict the order is arbitrary.)

\begin{cor}\label{cor:class}
Let $K$ be some Legendrian knot $L^\ell$, presented by the unknot stabilized $p$-times positively and $n$-times negatively, and with the rotation numbers on the two singular fibers taking the values
\[
\begin{array}{l}
\rot_1^i = \tb_1^i +1 \text{ for } i=0,\dots,K \text{ and } \rot_1^{K+1} \neq \tb_1^{K+1} +1 \\
\rot_2^i = -\tb_2^i -1 \text{ for } i=0,\dots,J \text{ and } \rot_2^{J+1} \neq -\tb_2^{J+1} -1.
\end{array}
\]
Denote by $D$ the denominator of $[a_1^0,\dots,a_1^{K-1}]^{-1}$ and by $D'$ the denominator of $[a_1^0,\dots,a_1^{K}]^{-1}$, and similarly, by  $E$ the denominator of $[a_2^0,\dots,a_2^{J-1}]^{-1}$ and by $E'$ the denominator of $[a_2^0,\dots,a_2^{J}]^{-1}$. (Whenever the continued fraction is empty, we formally take $[\,]=\frac{1}{0}$.)

\begin{itemize}[leftmargin=.5cm]
\item We have another presentation of the equivalent knot if either 
\[ D \leq p \text{ and there exists }  \begin{array}{l} j\leq J+1 \text{ and } p'\leq p_2^j: \\ \relax [a_1^0,\dots,a_1^{K-1}]^{-1}+[a_2^0,\dots,a_2^{j-1}, p']^{-1}=1, \end{array} \] 
or
 \[ E \leq n \text{ and there exists } \begin{array}{l} k\leq K+1 \text{ and } n'\leq n_1^k:\\ \relax [a_2^0,\dots,a_2^{J-1}]^{-1}+[a_1^0,\dots,a_1^{k-1},n']^{-1}=1.\end{array} \] 
In the first case it is given by
\[ 
\begin{array}{l} \ROT_1^{\,i}=-\tb_1^i-1 \text{ for } i=0,\dots,K \text{ and } \ROT_1^{\,K+1}=\rot_1^{K+1}-2\\
\ROT_2^{\,i}=\tb_2^i+1 \text{ for } i=0,\dots,j-1 \text{ and } \ROT_2^{\,j}=\rot_2^{j}-2p'\end{array} 
\]
and the rotation number of the knot decreased by $2D$,\\
in the second case by 
\[ 
\begin{array}{l} \ROT_2^{\,i}=\tb_1^i+1 \text{ for } i=0,\dots,J \text{ and } \ROT_2^{\,J+1}=\rot_2^{J+1}+2\\
\ROT_1^{\,i}=-\tb_2^i-1 \text{ for } i=0,\dots,k-1 \text{ and } \ROT_2^{\,k}=\rot_2^{k}+2n'\end{array} 
\]
and the rotation number of the knot increased by $2E$.

Any two presentations of the same knot are connected by a sequence of the described equivalences.

\item The knot $K$ is loose, if and only if one (and all) of its equivalent presentations satisfies 
\[ [a_1^0,\dots,a_1^{K}]^{-1}+[a_2^0,\dots,a_2^{J}]^{-1}<1 \]
and either \[ D'<E' \text{ and } p\geq D'\]
or \[E'<D' \text{ and } n\geq E'\]
(In particular, if no leading unknot of the singular fibers nor the knot itself is fully positive or fully negative, then $K$ is loose.)
\end{itemize}
\end{cor}

\begin{proof}
We rewrite relations from \cite[Section 5]{M.t} in accordance with the correspondence given in Theorem \ref{theorem:class}. Indeed, the equivalent presentations exactly correspond to isotopies (I1) from \cite[Proposition 5.2]{M.t}, and the loose knots to overtwisted surgeries given by (O1) in \cite[Proposition 5.1]{M.t}.

The simplification that looseness can be read from any given presentation is a consequence of the fact that the inequality $[a_1^0,\dots,a_1^{K}]^{-1}+[a_2^0,\dots,a_2^{J}]^{-1}<1$ fails only when presentation contains a balanced sublink, and in all other cases $D'$ and $E'$ are different.
\end{proof}

\begin{remark}[Negative torus knots in the standard structure]\label{rmk:Lstd}
In Proposition \ref{prop:tightC} we tell when the ambient contact structure of a Legendrian knots $L$, and hence of $L^\ell$, is tight; that is, if and only if the surgery diagram (of Figure \ref{fig:LSdgm}) contains a balanced sublink. From the classification in Theorem \ref{theorem:class} it follows that we can present in this form all Legendrian negative torus knots lying in $(S^3,\xi_\text{std})$, as classified by Etnyre and Honda \cite{EH} (recall from \cite{EH} that $\tb\leq -pq$).

Indeed, these are the only presentations for which the inequality in the second item of Corollary \ref{cor:class} is never satisfied. 
\end{remark}

\begin{remark}[Tension]\label{rmk:tension}
On the other hand, we know that a knot in an overtwisted structure becomes loose after sufficiently many stabilizations. However, we can read from the looseness criterion of Corollary \ref{cor:class} that for any non-loose realizations of $T_{(p,-q)}$ stabilizations of one sign suffice, while stabilizations of the opposite sign will keep the knot non-loose. Concretely, invoking the notation of Corollary \ref{cor:class}, the knot $K$ with $D'<E'$ has the positive tension equal to $t_+(K)=D'$ and the infinite negative tension $t_-(K)=\infty$, and oppositely,  the knot $K$ with $E'<D'$ has the negative tension equal to $t_-(K)=E'$ and the infinite positive tension $t_+(K)=\infty$.
\end{remark}

\begin{cor}[Coarse classification of transverse $T_{(p,-q)}$] \label{cor:transverse}
Any transverse realization $T$ of the knot $T_{(p,-q)}$ with Giroux torsion equal to zero arises as a transverse approximation of an $L^\ell$ for which $p<D'<E'$ in the notation of Corollary \ref{cor:class}.
\end{cor}

\begin{proof}
Recall that transverse knots can be thought of as Legendrian knots up to negative stabilizations. Hence, for a transverse knot to have tight complement, all negative stabilizations of its Legendrian push-off have to have tight complement as well. Since by stabilizing, the Thurston-Bennequin invariant eventually gets smaller than $-pq$, every transverse knot will be an approximation of some $L^\ell$. 
According to the preceding Remark \ref{rmk:tension}, the realizations with infinite negative tension has $D'<E'$, while $p<D'$ ensures that the Legendrian knot is non-loose.
\end{proof}

Note though that not all transversely non-loose $T_{(p,-q)}$ can arise as transverse approximation of an (non-stabilized) $L$ from Figure \ref{fig:LSdgm}. This can be observed from the isotopies in the first item of Corollary \ref{cor:class}: when $D$ is greater than 1, the $D$-times positively stabilized $L$ in the first presentation will be isotopic to $D$-times negatively stabilized knot in the second presentation, however the $d$-times positive stabilizations of $L$ for $d<D$ do not negatively destabilize to any non-loose knot with $\tb=-pq$ .

\bigskip
In the following, we will examine the question whether all elements of the minus version of the knot Floer homology can be realized as the Legendrian invariants (in the sense of Lisca, Ozsv\'ath, Stipsicz and Szab\'o \cite{LOSSz}) of a Legendrian realization in some contact structure on $S^3$ of the underlying knot. Notice that although the Legendrian invariant was introduced as a Legendrian isotopy invariant, its vanishing depends only on the contactomorphism type of the knot complement as can be understood from its sutured reinterpretation (established by Etnyre, Vela-Vick and Zarev in \cite{EVVZ}).

\begin{thm}\label{theorem:invL}
A Legendrian realization $K$ of the knot $T_{(p,-q)}$ in an overtwisted $S^3$ has non-zero Legendrian invariant $\mathfrak L(K)$ in $\HFK^-(T_{(p,q)})$, if and only if $K$ is up to negative stabilizations equivalent to a Legendrian knot $L^\ell$ with $p<D'<E'$. 

Hence, a transverse realization $T$ of the knot $T_{(p,-q)}$ in an overtwisted $S^3$ is strongly non-loose, if and only if it has non-zero transverse invariant $\mathfrak T(T)$ in $\HFK^-(T_{(p,q)})$. Moreover, they are classified by their transverse invariants, and actually, even by the self-linking number.
\end{thm}

\begin{proof}
First, if a knot $K$ has a non-zero invariant $\mathfrak L(K)$, so do all its negative stabilizations. Hence, a Legendrian negative torus knot has non-zero invariant only if it is strongly non-loose as a transverse knot. Indeed, in the presence of Giroux $2\pi$-torsion all Heegaard Floer invariants vanish \cite{SV}, while the vanishing of invariants when we add boundary parallel $\pi$-torsion is a consequence of an odd change (by $\slk=2A-1$) in $d_3$ and the structure of $\HFK^-(T_{(p,q)})$ (see Remark \ref{rmk:d3}). So, by Corollary \ref{cor:transverse}, $K$ is up to negative stabilizations equivalent to a Legendrian knot $L^\ell$ with $p<D'<E'$. 

So, we just need to prove that every $L^\ell$ with $p<D'<E'$ (hence, every transversely non-loose realization) has non-zero invariant, and that they distinguish non-equivalent realizations.

Nonvanishing is a direct application of \cite[Theorem 1.1]{M.s}. Indeed, the negative tension of $L^\ell$ with $p<D'<E'$ is infinite, and by the correspondence of Theorem \ref{theorem:class} the $-n$-contact surgeries of all stabilizations negative along such $L^\ell$ are tight for all $n\in\NN$. Since the underlying smooth manifolds are small Seifert fibered $L$-spaces with Euler number $-1$, the tight structures also have non-zero contact invariant (by \cite[Corollary 1.4]{M.t}). Hence, the theorem applies and $\mathfrak L(L^\ell)\neq0$, and then $\mathfrak L(K)\neq0$ for any Legendrian $K$ which is equivalent to $L^\ell$ up to negative stabilizations.

Nonequivalent $L^\ell$ give by the correspondence of Theorem \ref{theorem:class} non-isotopic tight structures on $M(-1;\frac{p-p'}{p},\frac{q'}{q},\frac{2}{2\ell+1})$; by \cite[Theorem 1.3]{M.t} also their induced $\spinc$-structures are different, and so are the relative $\spinc$-structures of the knots. Therefore, nonequivalent transverse knots are differed (already) by the Alexander grading of their transverse invariants, or equivalently, by their self-linking number.
\end{proof}

\begin{remark}\label{rmk:o=t}
Applying Remark \ref{rmk:tension} we see that any Legendrian realization $K$ of the knot $T_{(p,-q)}$ in an overtwisted structure on $S^3$ has the order equal to the tension (in terminology of Baker and Onaran \cite{BO}).
\end{remark}

\begin{remark}[Ambient overtwisted structures] \label{rmk:d3}
Here we show in which overtwisted structures transverse strongly non-loose negative torus knots $T_{(p,-q)}$ appear, and that each of these structures can be occupied by only one such knot up to stabilizations.

Nonvanishing of Legendrian invariants makes it easy to see that transverse strongly non-loose negative torus knots $T_{(p,-q)}$, and in fact all $L$ from Figure \ref{fig:LSdgm}, live in the overtwisted structures with the $3$-dimensional invariant positive, even, and at most $(p-1)(q-1)$. Indeed, the knot Floer homology of a positive torus knot \cite{OSz.L} takes the following form 
\[ \HFK^-(T_{(p,q)})\cong \FF[U]_{(A_0,M_0)} \oplus \bigoplus_i \big(\FF[U]/(U^{n_i})\big)_{(A_i,M_i)} \]
where $(A_i,M_i)$ is the bigrading of the generator for each summand. Now, since $\HFK^-(T_{(p,q)})$ is the homology of the associated graded complex of an upper staircase, the generators of the direct summands correspond exactly to the bottom dots of the staircase and $2A_i-M_i$ (which is preserved by $U$-multiples, and hence throughout each summand) is equal to twice the second coordinate of the bottom dots of the upper staircase (when based on the axes of the first quadrant). Therefore,  $\HFK^-(T_{(p,q)})$ is supported only in finitely many values of $2A-M$, all of which are even and in the interval $[0, (p-1)(q-1)]$. This finishes the proof for the transversely non-loose knots whose invariants are non-zero by Theorem \ref{theorem:invL}; however, from Remark \ref{rmk:tension} we see that every $L$ is either transversely non-loose or it gets such after all unknots in its surgery presentation are taken with its reversed orientation, which preserves the $3$-dimensional invariant.

Moreover, since by Theorem \ref{theorem:invL} transversely non-loose $T_{(p,-q)}$ are classified by the transverse invariants, the only transverse knots living in the same $3$-dimensional invariant are the ones related by (positive) stabilizations, which (for their Legendrian approximations) are exactly described by the first item of Corollary \ref{cor:class}. This is because the equality $2A_i - M_i = 2A_j - M_j$ holds only if $i$ and $j$ are equal, and so, the same $d_3$ is shared only by the knots whose transverse invariants are in the same summand. But since we already know that different transverse $T_{(p,-q)}$ cannot have the same invariant, all the ones with the invariants in the same summand are related by stabilizations.

In particular, we notice that the single generator in the top Alexander grading $A=\frac{(p-1)(q-1)}{2}$ is always realized by the binding for the open book supporting the overtwisted structure with $d_3=2A$ (as proved by Vela-Vick in \cite{VV}), and this is the only transversely non-loose $T_{(p,-q)}$ living in the maximal $d_3$. In our presentation of Figure \ref{fig:LSdgm} the corresponding knots are given by all the surgery curves fully positive, that is, when $\rot_i^j=-\tb_i^j-1$ for all the surgery curves. (Indeed, the maximal $d_3$ is attained at the maximal $\rot$-vector because the coefficients in the inverse of the intersection matrix are all non-negative, by an analogous argument to \cite[Section 4.1]{LiSi}.)
\end{remark}

\begin{example}\label{example}
We demonstrate that not all the torsion elements of $\HFK^-(T_{(p,q)})$ are realized as Legendrian invariants of (transversely) non-loose negative torus knots $T_{(p,-q)}$. The first example is a (very special) family for which this is the case, the second is a family for which all $1$-torsion elements are realized, and the third is a concrete example for which not all $1$-torsion elements are realized, but there is a Legendrian invariant of order $2$.

\begin{itemize}[leftmargin=.5cm]

\item $T_{(2,-2n+1)}$: The two singular fibers have coefficients $\frac{1}{2}=[2]^{-1}$ and $\frac{n}{2n-1}=[2,n]^{-1}$. There are $n-1$ strongly non-loose representatives, distinguished by the rotation number on the $(-n)$-framed surgery. Meanwhile, the knot Floer homology takes the form $\HFK^-(T_{(2,2n-1)})\cong \FF[U]\oplus\FF^{n-1}$. Also, the set of pairs of Alexander and Maslov gradings agrees with the set of triples $(\tb,\rot,d_3)$ for the listed knots. This has already been observed in \cite[Remark 6.11]{LOSSz}: in particular, our Theorem \ref{theorem:invL} confirms the conjecture of Lisca, Ozsv\'ath, Stipsicz and Szab\'o that the Legendrian invariants of the knots $L_{k,l}$ (in the notation of \cite{LOSSz}) are non-zero, and hence present generators for all $U$-torsion elements of $\HFK^-(T_{(2,2n-1)})$.

\item $T_{(n,-n-1)}$: The two singular fibers have coefficients $\frac{1}{n}=[n]^{-1}$ and $\frac{n}{n+1}=[2^{\times n}]^{-1}$. The number of the relevant non-loose representatives (with $D'<E'$, in the notation of Corollary \ref{cor:class}) is $n-1$ (the rotation number on the $(-3)$-framed unknot is 1, and on the $(-n-1)$-framed unknot it is any of $\{n-1,n-3,\dots,-n+3\}$), while the knot Floer homology takes the form \[\HFK^-(T_{(n,n+1)})\cong \FF[U]\oplus \bigoplus_{i=1}^{n-1}\FF[U]/(U^i).\] It is possible to check that the bigradings computed from $(\tb,\rot,d_3)$ agree with the bigradings of the bottoms of the torsion summands.

\item $T_{(5,-8)}$: The two singular fibers have coefficients $\frac{5}{8}=[2,3,2]^{-1}$ and $\frac{2}{5}=[3,2]^{-1}$. The number of the relevant non-loose representatives (with $D'<E'$, in the notation of Corollary \ref{cor:class}) is $5$. The knot Floer homology takes the form 
\[ \HFK^-(T_{(5,8)})\cong\FF[U]\oplus\FF[U]/(U^4)\oplus \left(\FF[U]/(U^2)\right)^2 \oplus \left(\FF[U]/(U)\right)^6, \]
the Legendrian invariants lie at $(A,M)$ equal to:
\[\begin{array}{l}
(14,0) \text{ in }  \left(\FF[U]/(U)\right)_{(14,0)}, \text{ for } \rot=[1,1,0; 2,0]\\
(4,-6) \text{ in }  \left(\FF[U]/(U)\right)_{(4,-6)}, \text{ for } \rot=[1,1,0; 0,0]\\
(-2,-12) \text{ in } \left(\FF[U]/(U^2)\right)_{(-1,-10)},\text{ for } \rot=[1,-1,0; 2,0]\\
(-12,-26), (-11,-24) \text{ in }  \left(\FF[U]/(U^4)\right)_{(-9,-20)},\\
\ \ \ \ \ \ \ \ \ \ \ \ \ \ \ \ \ \ \ \ \ \ \ \ \ \ \ \ \ \ \ \ \ \ \ \ \ \ \ \ \text{ for } \rot=[1,-1,0;0,0], [-1,1,0; 2,0]
\end{array}\]
Note that the last two presentations are related by isotopy of Corollary \ref{cor:class} when the knot is stabilized once negatively, respectively once positively.
\end{itemize} 
\end{example}

\bigskip
Finally, we look back at the tight structures on \[ S^3_{-pq-m}(T_{(p,-q)})=M(-1;\frac{p-p'}{p},\frac{q'}{q},\frac{1}{m+1}) \] Here, we will notice that the number of tight structures on these manifolds stabilizes once $m$ gets big enough. Precisely, for $m>q$ the manifold $S^3_{-pq-m-1}(T_{(p,-q)})$ admits one more tight contact structure than the manifold $S^3_{-pq-m}(T_{(p,-q)})$. 

Indeed, recall first that the fillable structures always contain a balanced sublink \cite{M.f}. Thus, new relations among seemingly different contact presentations (of fillable structures) are induced by isotopies (I1) of \cite[Proposition 5.2]{M.t}, and they are newly appearing only as long as the denominators of the first two Seifert constants are bigger than $m$. Once $m>q$, all equivalences are already established and the one more structure is always coming from the choice of the one more stabilization on the knot supporting the third singular fiber. 
Thinking about Legendrian representations of $T_{(p,-q)}$ instead, we have  observed already in Remark \ref{rmk:Lstd} that $L^m$ given by the contact surgery presentations containing a balanced sublink correspond exactly to Legendrian realizations in the standard contact $S^3$. The equivalences between the knot presentations from Corollary \ref{cor:class} are exactly the isotopies (I1) mentioned above for the fillable structures. In fact, we can see directly from the mountain range that for the low enough Thurston-Bennequin invariants, the number of Legendrian realizations of $T_{(p,-q)}$ in the standard contact structure on $S^3$ increases by one as the Thurston-Bennequin number decreases by one. This is indeed the case for $\tb\leq -pq-q$.  Therefore, for $m>q$ we established a one to one correspondence between fillable structures on $S^3_{-pq-m}(T_{(p,-q)})$ and transverse realizations in the standard structure with the self-linking greater than or equal to $-pq-q+p-m$.

On the other hand, we have the following relation between the non-loose realizations of the knot  $T_{(p,-q)}$ in overtwisted structures on $S^3$ and the non-fillable tight contact structures on the large negative surgeries.

\begin{prop}\label{prop:number}
The number of non-fillable tight contact structures on the $(-m)$-surgery along $T_{(p,-q)}$ for large $m\in\NN$ equals 
twice the number of the strongly non-loose transverse realizations of the knot $T_{(p,-q)}$ decreased by the number of overtwisted structures in which they appear.  This is further equal to twice the number of torsion elements in $\HFK^-(T_{(p,q)})$  which are realized as the Legendrian invariants, minus the number of realized $1$-torsion elements.
\end{prop}

\begin{proof}
We have observed in Theorem \ref{theorem:class} that the Legendrian realizations without Giroux torsion of $T_{(p,-q)}$ with $\tb=-pq-m$  are in one to one correspondence with tight structures on $M(-1;\frac{p-p'}{p},\frac{q'}{q},\frac{2}{2m+1})$, and that taking a single Legendrian surgery along them gives rise to all tight structure on $S^3_{-pq-m-1}(T_{(p,-q)})$. Above we have already related the fillable structures with the knots in the standard $S^3$. So, it only remains to prove when the Legendrian surgery along different non-loose Legendrian knots results in the same (non-fillable) tight structure.

Comparing the equivalences from Corollary \ref{cor:class} with the isotopies described in \cite[Section 5]{M.t}, we see that the only identifications of tight structures (not arising from equivalences of knots) are provided by  isotopies (I2) of \cite[Proposition 5.2]{M.t}. These isotopies connect surgery presentations with either a fully negative or a fully positive unknot representing the third singular fiber and neither  leading unknot of the other two singular fibers fully negative or positive, respectively. They identify them in pairs: the one with the fully negative unknot on the third fiber and the one with this unknot fully positive which has rotation numbers of the other leading unknots decreased by two and the rotation numbers of all the other surgery curves the same. According to  Remark \ref{rmk:tension}, the knots supporting the third singular fiber in these presentations are exactly the non-loose Legendrian $T_{(p,-q)}$ which have the positive, respectively negative, tension equal to one. Furthermore, according to Remark \ref{rmk:o=t} also their order equals one. Therefore, we need to subtract from the number of Legendrian realizations (which is twice the number of transverse realizations, and equals the number of realized torsion $\HFK$-classes)  the number of $1$-torsion elements, realized by Legendrian invariants, or equally, the number of overtwisted structures that admit a non-loose realization.
\end{proof}

%%%%%%%%%%%%%%%%%%%%%%%%%%%%%%%%%%%%%%%%%%%%%%%%%%%%%%%%%%%%%%

%%%%%%%%%%%%%%%%%%%%%%%%%%%%%%%%%%%%%%%%%%%%%%%%%%%%%%%%%%%%%%

%%%%%%%%%%%%%%%%%%%%%%%%%%%%%%%%%%%%%%%%%%%%%%%%%%%%%%%%%%%%%%
%%%%%%%%%%%%%%%%%%%%%%%%%%%%%%%%%%%%%%%%%%%%%%%%%%%%%%%%%%%%%%
\end{document}